\documentclass{amsart}

\usepackage{amssymb}

\newtheorem{thm}{Theorem}

\newtheorem{question}[thm]{Question}
\newtheorem{lem}[thm]{Lemma}
\newtheorem{cor}[thm]{Corollary}

\newtheorem{prop}[thm]{Proposition}

\theoremstyle{definition}
\newtheorem{defn}[thm]{Definition}

\newtheorem{say}[thm]{}
\newtheorem{exmp}[thm]{Example}

\newtheorem{rem}[thm]{Remark}          

\newtheorem{ack}{Acknowledgments}

\newtheorem{defn-thm}[thm]{Definition--Theorem}  
\newtheorem{defn-lem}[thm]{Definition--Lemma}  

\theoremstyle{remark}

\renewcommand{\c}[0]{{\mathbb C}}  

\renewcommand{\o}[0]{{\mathcal O}} 
\newcommand{\z}[0]{{\mathbb Z}}

\renewcommand{\a}[0]{{\mathbb A}}

\newcommand{\p}[0]{{\mathbb P}}
\newcommand{\f}[0]{{\mathbb F}}

\newcommand{\qtq}[1]{\quad\mbox{#1}\quad}
\newcommand{\spec}[0]{\operatorname{Spec}}

\newcommand{\red}[0]{\operatorname{red}}    
    
\newcommand{\im}[0]{\operatorname{im}}

\newcommand{\chr}[0]{\operatorname{char}}

\newcommand{\mor}[0]{\operatorname{Mor}}

\def\into{\DOTSB\lhook\joinrel\to}

\begin{document}
\bibliographystyle{amsalpha}

\title{Quotients by finite equivalence relations}
\author{J\'anos Koll\'ar\\{\ }\\
{\scriptsize with an appendix by} C.\ Raicu
}


\maketitle

Let $f:X\to Y$ be a finite  morphism of schemes.
Given $Y$, one can easily describe $X$ by the coherent sheaf
of algebras $f_*\o_X$. 
Here our main interest is the converse.
Given $X$, what kind of data do we need to construct $Y$?
For this question,  the surjectivity of $f$ is indispensable.

The fiber product $X\times_YX\subset X\times X$
defines an equivalence relation on $X$, and one might hope
 to reconstruct $Y$ as the quotient of $X$ by this
equivalence relation. Our main interest is in the cases when
$f$ is not flat. 
A typical example we have in mind is when   $Y$ is not normal and $X$ is
its normalization. In these cases,
the fiber product $X\times_YX$ can be  rather complicated.
Even if $Y$ and $X$ are  pure dimensional and CM, 
$X\times_YX$ can have irreducible components of different
dimension and its connected components need not be
pure dimensional.
None of these difficulties appear if $f$ is flat 
\cite{MR0232781, sga3} or
if $Y$ is normal (\ref{quot.pure.dim.lem}).

The aim of this note is to give many examples, review known results,
pose questions and to prove a few  theorems
concerning finite equivalence relations.

\section{Definition of equivalence relations}

\begin{defn}[Equivalence relations]\label{eq.rel.defn}
 Let $X$ be an $S$-scheme
and  $\sigma:R\to X\times_SX$  a morphism
(or $\sigma_1,\sigma_2:R\rightrightarrows X$ a pair of morphisms).
We say that  $R$ is an {\it equivalence relation}
on $X$ if, for every scheme $T\to S$,
we get a (set theoretic)  equivalence relation 
$$
\sigma(T):\mor_S(T,R)\into
\mor_S(T,X)\times \mor_S(T,X).
$$
Equivalently, the following conditions hold:
\begin{enumerate}
 \item $\sigma$ is a monomorphism (\ref{monom.defn})
\item (reflexive) $R$ contains the diagonal $\Delta_X$.
\item  (symmetric) There is an involution $\tau_R$ on $R$ such that
$\tau_{X\times X}\circ\sigma\circ\tau_R=\sigma$,
where $\tau_{X\times X}$ denotes the involution which interchanges the
two factors of $X\times X$. 
\item (transitive) For $1\leq i<j\leq 3$ set $X_i:=X$ and let
$R_{ij}:=R$ when it maps to $X_i\times_SX_j$.
Then  the coordinate projection 
of $R_{12}\times_{X_2}R_{23}$ to $X_1\times_SX_3$
factors through $R_{13}$:
$$
R_{12}\times_{X_2}R_{23}\to R_{13}\stackrel{\pi_{13}}{\longrightarrow}
X_1\times_SX_3.
$$
\end{enumerate}

We say that  $\sigma_1,\sigma_2:R\rightrightarrows X$
is a {\it finite} equivalence relation
if the maps $\sigma_1,\sigma_2$ are finite. In  this case,
$\sigma:R\to X\times_SX$ is also finite, hence a closed embedding
(\ref{monom.defn}).
\end{defn}

\begin{defn}[Set theoretic equivalence relations]\label{setth.eq.rel.defn}
Let $X$ and $R$ be  reduced $S$-schemes. We say that  a morphism
 $\sigma:R\to X\times_SX$ 
  is a {\it set theoretic equivalence relation}
on $X$ if, for every geometric point $\spec K\to S$,
we get an equivalence relation on $K$-points
$$
\sigma(K):\mor_S(\spec K,R)\into
\mor_S(\spec K,X)\times \mor_S(\spec K,X).
$$
Equivalently, 
\begin{enumerate}
 \item $\sigma$ is geometrically injective.
\item (reflexive) $R$ contains the  diagonal $\Delta_X$.
\item  (symmetric) There is an involution $\tau_R$ on $ R$ such that
$\tau_{X\times X}\circ\sigma\circ\tau_R=\sigma$
where $\tau_{X\times X}$ denotes the involution which interchanges the
two factors of $X\times X$.
\item (transitive) For $1\leq i<j\leq 3$ set $X_i:=X$ and let
$R_{ij}:=R$ when it maps to $X_i\times_SX_j$.
Then  the coordinate projection 
of $\red\bigl(R_{12}\times_{X_2}R_{23}\bigr)$ to $X_1\times_SX_3$
factors through $R_{13}$:
$$
\red\bigl(R_{12}\times_{X_2}R_{23}\bigr)\to R_{13}
\stackrel{\pi_{13}}{\longrightarrow}
X_1\times_SX_3.
$$
\end{enumerate}
Note that the  fiber product need not be reduced,
and taking the reduced structure above is essential,
as shown by (\ref{nonred.noneq.exmp}).

 It is sometimes convenient to consider  finite  morphisms
 $p:R\to X\times_SX$
such that the injection $i:p(R)\into  X\times_SX$
is a set theoretic equivalence relation.
Such a $p:R\to X\times_SX$ is called a 
   {\it set theoretic pre-equivalence relation.}
\end{defn}

\begin{exmp} \label{nonred.noneq.exmp} On $X:=\c^2$ consider the $\z/2$-action
$(x,y)\mapsto (-x,-y)$. This can be given by a 
set theoretic equivalence relation
$R\subset X_{x_1,y_1}\times X_{x_2,y_2}$
defined by the ideal
$$
(x_1-x_2,y_1-y_2)\cap (x_1+x_2,y_1+y_2)=
(x_1^2-x_2^2, y_1^2-y_2^2, x_1y_1-x_2y_2,  x_1y_2-x_2y_1)
$$
in $\c[x_1,y_1,x_2,y_2]$.
We claim that this is {\it not} an equivalence relation.
The problem is with transitivity.
The defining ideal of  $R_{12}\times_{X_2}R_{23}$ in
$\c[x_1,y_1,x_2,y_2,x_3,y_3]$ is
$$
(x_1^2-x_2^2, y_1^2-y_2^2, x_1y_1-x_2y_2,  x_1y_2-x_2y_1,
x_2^2-x_3^2, y_2^2-y_3^2, x_2y_2-x_3y_3,  x_2y_3-x_3y_2).
$$
This contains $(x_1^2-x_3^2, y_1^2-y_3^2, x_1y_1-x_3y_3)$
but it does not contain $ x_1y_3-x_3y_1$.
Thus there is no map $R_{12}\times_{X_2}R_{23}\to R_{13}$.
Not, however, that the problem is easy to remedy.
Let $R^*\subset X\times X$ be
defined by the ideal
$$
(x_1^2-x_2^2, y_1^2-y_2^2, x_1y_1-x_2y_2)\subset \c[x_1,y_1,x_2,y_2].
$$
We see that $R^*$ defines an equivalence relation.
The difference between $R$ and $R^*$ is one embedded point at the
origin.
\end{exmp}

\begin{defn}[Categorical and geometric quotients]\label{cat.geom.quot}
Given two morphisms, $\sigma_1,\sigma_2:R\rightrightarrows X$,
 there is at most one  scheme $q:X\to (X/R)^{cat}$ such that
$q\circ \sigma_1=q\circ \sigma_2$ and $q$ is
universal with this property.
 We call $(X/R)^{cat}$ the 
{\it categorical quotient} (or {\it coequalizer}) of 
$\sigma_1,\sigma_2:R\rightrightarrows X$.

The categorical quotient is  easy to construct in the affine case.
Given $\sigma_1,\sigma_2:R\rightrightarrows X$, the categorical quotient
$(X/R)^{cat}$
 is the spectrum of the
$S$-algebra
$$
\ker\Bigl[\o_X\stackrel{\sigma_1^*-\sigma_2^*}{\longrightarrow} \o_R\Bigr].
$$

Let $\sigma_1,\sigma_2:R\rightrightarrows X$  be 
 a finite  equivalence relation. 
We say that   $q:X\to Y$ is a
 {\it geometric quotient} of $X$ by $R$
if  
\begin{enumerate}
\item   $q:X\to Y$ is the categorical quotient $q:X\to (X/R)^{cat}$,
\item  $q:X\to Y$ is finite, and
\item for every geometric point $\spec K\to S$,
the fibers of 
$q_K:X_K(K)\to Y_K(K)$
are the
$\sigma\bigl(R_K(K)\bigr)$-equivalence classes of 
 $X_K(K)$.
\end{enumerate}

The geometric quotient is denoted by $X/R$.
\end{defn}

The main example to keep in mind is the following,
which easily follows from (\ref{quot.X/S.finite.lem}) 
and the construction of
$(X/R)^{cat}$ for affine schemes.

\begin{exmp}\label{basic.quot.exmp}
 Let $f:X\to Y$ be a finite and surjective
morphism. Set $R:=\red(X\times_YX)\subset X\times X$
and let $\sigma_i:R\to X$ denote the coordinate projections.
Then the geometric quotient $X/R$ exists and
$X/R\to Y$ is a finite and universal homeomorphism (\ref{univ.homeo.defn}).
Therefore, if $X$ is the 
normalization of $Y$, then  $X/R$ is the weak normalization of $Y$.
(See \cite[Sec.7.2]{rc-book} for basic results on semi-normal 
and weakly normal schemes.)

By taking  the reduced structure of
$X\times_YX$ above, we chose to  focus on the 
set-theoretic properties of $Y$. 
However, as (\ref{sch.quot.exmp}) shows, even if
$X,Y$ and $X\times_YX$ are all reduced,
$X/R\to Y$ need not be an isomorphism.
Thus $X$ and $X\times_YX$ do not determine $Y$ uniquely.
\end{exmp}

In Section \ref{first.exmp.sec} we give examples  of  finite,
set theoretic equivalence relations $R\rightrightarrows X$
such that the categorical quotient $(X/R)^{cat}$ is non-Noetherian
and there is no geometric quotient. This can happen even when $X$ is very nice,
for instance a smooth variety over $\c$.
Some elementary results about the existence of
geometric quotients are discussed in Section \ref{basic.res.sec}.

An inductive plan to construct geometric quotients is outlined in
Section \ref{induct.plan.sect}.
As an application, we prove in Section \ref{pos.char.sec}
the following:

\begin{thm}\label{quot.by.R.charp} Let $S$ be a Noetherian $\f_p$-scheme and
$X$ an algebraic space which is essentially of finite type over $S$.
Let $R\rightrightarrows X$ be a finite, set theoretic equivalence relation.
Then the geometric  quotient $X/R$ exists.
\end{thm}

\begin{rem} There are many algebraic spaces which are not of finite
type and such  that the Frobenius map
$F^q:X\to X^{(q)}$ is  finite. By a result of Kunz (see
\cite[p.302]{mats-ca}) such algebraic spaces are excellent. 
As the proof  shows, 
(\ref{quot.by.R.charp})  remains valid for algebraic spaces
satisfying this property. 
\end{rem}

In the Appendix, C.~Raicu constructs finite scheme theoretic
equivalence relations $R$ on $X=\a^2$  (in any characteristic)
such that the geometric quotient $X/R$ exists
yet $R$ is strictly smaller than the fiber
product $X\times_{X/R}X$.
Closely related examples are in \cite{venken, philippe}.

In characteristic zero, this leaves open the following:

\begin{question}\label{sch.th.quot.quest}
 Let $R\subset  X\times X$ be a 
scheme theoretic equivalence relation
such that the coordinate projections  $R\rightrightarrows X$ are finite.

Is there a geometric quotient $X/R$?
\end{question}

A special case of the quotient problem, called gluing or pinching,
is discussed in Section \ref{glue.sec}.
This follows the works of \cite{artin70},
\cite{ferrand} (which is based on an unpublished manuscript from '70)
and \cite{raoult}.

\section{First examples}\label{first.exmp.sec}

The next examples show that in many cases, the
categorical quotient of a very nice scheme $X$ can be non-Noetherian.
We start with a nonreduced example and then
we build it up to smooth ones.

\begin{exmp} \label{first.exmp}
Let $k$ be a field and consider $g_i:k[x,\epsilon]\to k[x,\epsilon]$ 
where 
$$
g_1\bigl(a(x)+\epsilon b(x)\bigr)=a(x)+\epsilon b(x)
\qtq{and}
g_2\bigl(a(x)+\epsilon b(x)\bigr)=a(x)+\epsilon \bigl(b(x)+a'(x)\bigr).
$$

If $\chr k=0$ then the coequalizer is the spectrum of
$$
\ker\Bigl[k[x,\epsilon]\stackrel{g_1^*-g_2^*}{\longrightarrow} 
k[x,\epsilon]\Bigr]=k+\epsilon k[x].
$$
Note that $k+\epsilon k[x]$ is not Noetherian
and its only prime ideal is $\epsilon k[x]$.

If $\chr k=p$ then the coequalizer
is the spectrum of the finitely generated $k$-algebra
$$
\ker\Bigl[k[x,\epsilon]\stackrel{g_1^*-g_2^*}{\longrightarrow} 
k[x,\epsilon]\Bigr]=k[x^p]+\epsilon k[x].
$$
\end{exmp}

It is not surprising that set theoretic equivalence relations
behave badly on nonreduced schemes. However, the above example
is easy to realize on reduced and even on smooth schemes.

\begin{exmp} (cf.\ \cite[p.342]{holmann}) 
Let $p_i:Z\to Y_i$ be finite morphisms for $i=1,2$.
 We can construct out of them  an
equivalence relation on $Y_1\amalg Y_2$
where $R$ is the union of the diagonal with
two copies of $Z$, one of which maps as 
$$
(p_1,p_2):Z\to Y_1\times Y_2\subset \bigl( Y_1\amalg Y_2\bigr)
\times
\bigl( Y_1\amalg Y_2\bigr),
$$
the other its symmetric pair. The categorical quotient
$\bigl(\bigl( Y_1\amalg Y_2\bigr)/R\bigr)^{cat}$ is also 
the universal push-out of $Y_1\stackrel{p_1}{\leftarrow} Z
  \stackrel{p_2}{\to}Y_2$.
If $Z$ and the $Y_i$ are affine over $S$, then
it  is the spectrum of the
$S$-algebra
$$
\ker\Bigl[\o_{Y_1}+\o_{Y_2}\stackrel{p_1^*-p_2^*}{\longrightarrow} \o_Z\Bigr].
$$

For the first example 
let $Y_1\cong Y_2:=\spec  k[x,y^2,y^3]$
and $Z:=\spec k[u,v]$ with $p_i$ given by
$$
p_1^*: (x,y^2,y^3)\mapsto (u,v^2,v^3)
\qtq{and}
p_2^*:(x,y^2,y^3)\mapsto (u+v,v^2,v^3).
$$
Since the $p_i^*$ are injective,
 the categorical quotient is the spectrum of the
$k$-algebra
$k[u,v^2,v^3]\cap k[u+v,v^2,v^3]$. 
Note that
$$
\begin{array}{lll}
k[u,v^2,v^3] & = & \bigl\{f_0(u)+
\sum_{i\geq 2} v^i f_i(u)\ :\ f_i\in k[u]\bigr\}\qtq{and}\\
k[u+v,v^2,v^3] & = & \bigl\{f_0(u)+vf'_0(u)+\sum_{i\geq 2} v^i f_i(u)
\ :\ f_i\in k[u]\bigr\}.
\end{array}
$$
As in  (\ref{first.exmp}),
 if $\chr k=0$ then the categorical quotient is the spectrum of the
non-Noetherian algebra
$k+\sum_{n\geq 2}v^nk[u]$.
 If $\chr k=p$ then the geometric quotient is given by
the finitely generated $k$-algebra
$$
k[u^p]+\sum_{n\geq 2}v^nk[u].
$$
\end{exmp}

This example can be embedded into a
set theoretic equivalence relation on a smooth variety.

\begin{exmp}
Let $Y_1\cong Y_2:=\a^3_{xyz}$, $Z:=\a^2_{uv}$ and 
$$
p_1^*: (x_1,y_1,z_1)\mapsto (u,v^2,v^3)
\qtq{and}
p_2^*(x_2,y_2,z_2)\mapsto (u+v,v^2,v^3).
$$
By the previous computations, in
characteristic zero the categorical quotient
is given by
$$
k+(y_1,z_1)+(y_2, z_2)\subset k[x_1,y_1,z_1]+k[x_2,y_2,z_2],
$$
where $(y_i,z_i)$ denotes the ideal
$(y_i,z_i)\subset k[x_i,y_i,z_i]$. A minimal generating set is given by
$$
y_1x_1^m, z_1x_1^m, y_2x_2^m, z_2x_2^m\ :\ m=0,1,2, \dots
$$

In positive characteristic the
categorical quotient
is given by
$$
k[x_1^p,x_2^p]+ (y_1,z_1)+(y_2, z_2)\subset k[x_1,y_1,z_1]+k[x_2,y_2,z_2].
$$
A minimal generating set is given by
$$
x_1^p, x_2^p, y_1x_1^m, z_1x_1^m, y_2x_2^m, z_2x_2^m\ :\ m=0,1,\dots, p-1.
$$

\end{exmp}

\begin{exmp}\label{nagata.exmp}
The following example, based on  \cite{Nagata69}, shows that
even for rings of invariants of finite group actions
some finiteness assumption on $X$  is necessary
in order to obtain geometric quotients.

 Let $k$ be a field of characteristic
$p>0$  and $K:= k(x_1,x_2,\dots)$, where the
$x_i$ are algebraically independent over $k$.  Let 
$$
D:= \sum_ix_{i+1}\frac{\partial}{\partial x_i}\qtq{be a derivation of $K$.}
$$
Let $F:=\{f\in K\vert D(f)=0\}$ be the subfield of constants.
Set 
$$
R=K+\epsilon K\qtq{where $\epsilon^2=0$ and} \sigma: f+\epsilon g\mapsto
f+\epsilon (g+D(f)).
$$
$R$ is a local Artin ring.
It is easy to check that $\sigma$ is an automorphism of $R$ of order $p$. 
The fixed ring is $R^{\sigma}=F+\epsilon K$. Its maximal ideal is
$m:= (\epsilon K)$ and  generating sets of $m$ correspond to   $F$-vectorspace
bases
of $K$. Next we show that 
the $x_i$ are linearly independent
over $F$ which implies that   $R^{\sigma}$ is not
Noetherian.

 Assume that we have a relation $\sum_{i\leq n} f_ix_i=0$
We may assume that $f_n=1$ and
$f_i\in F\cap k(x_1,\dots,x_r)$ for some $r$.
Apply $D$ to
get that
$$
0=\sum_{i\leq n} f_iD(x_i)=\sum_{i\leq n} f_ix_{i+1}.
$$ 
Repeating $s$ times gives that 
$\sum_{i\leq n} f_ix_{i+s}=0$, or, equivalently
$$
x_{n+s}=-\sum_{i\leq n-1} f_ix_{i+s}.
$$
This is impossible
if $n+s>r$; a contradiction. 

It is easy to see that $R$ is not a submodule of any
finitely generated  $R^{\sigma}$-module.
\end{exmp}

\begin{exmp} This example of \cite{naga}
gives a 2-dimensional regular local ring $R$ and an
automorphism of order 2 such that the
ring of invariants is not Noetherian.

Let $k$ be a field of characteristic 2
 and $K:= k(x_1, y_1, x_2, y_2, \dots)$, where the
$x_i, y_i$ are algebraically independent over $k$.  Let 
$R:=K[[u,v]]$ be the power series ring in 2 variables.
Note that  $R$ is a 2-dimensional regular local ring, but it is not
 essentially of finite type over $k$.
Define a  derivation of $K$ to $R$ by
$$
D_K:= \sum_i v(x_{i+1}u+y_{i+1}v)\frac{\partial}{\partial x_i} +
u(x_{i+1}u+y_{i+1}v)\frac{\partial}{\partial y_i}.
$$
This extends to a  derivation of $R$ to $R$ by
setting $D_R|_K=D_K$ and $D_R(u)=D_R(v)=0$.
Note that $D_R\circ D_R=0$, thus
$\sigma: r\mapsto r+D_R(r)$ is an order 2 automorphism of $R$.
We claim that the ring of invariants $R^{\sigma}$ is not Noetherian.

To see this, note first that $x_{i}u+y_{i}v\in R^{\sigma}$
for every $i$.
\medskip

{\it Claim.}  For every $n$,
$x_{n+1}u+y_{n+1}v\not\in \bigl(x_{1}u+y_{1}v, \dots, x_{n}u+y_{n}v\bigr)
R^{\sigma}$.
\medskip

Proof. Assume the contrary and write
$$
x_{n+1}u+y_{n+1}v=\sum_{i\leq n} r_i\bigl(x_{i}u+y_{i}v\bigr)
\qtq{where $r_i\in R^{\sigma}$.}
$$
Working modulo $(u,v)^2$ and gathering the terms involving $u$,
we get an equality
$$
x_{n+1}\equiv \sum_{i\leq n}  r_ix_{i} \mod R^{\sigma}\cap (u,v)R.
$$
Applying $D_R$ and again gathering the terms involving $u$ we
obtain that
$x_{n+2}
\equiv\sum_{i\leq n}  r_ix_{i+1}$
modulo $R^{\sigma}\cap (u,v)R$.
Repeating this $s$ times gives
$$
x_{n+s+1}=\sum_{i\leq n}  \bar r_ix_{i+s}\qtq{where $\bar r_i\in K$.}
$$
Since the $\bar r_i$ involve only finitely many variables,
we get a contradiction for large $s$. Thus
$$
(p_1)\subset (p_1,p_2)\subset (p_1,p_2,p_3)\subset \cdots
$$
is an infinite increasing sequence of ideals in $R^{\sigma}$.
\end{exmp}

The next examples  show that,  if $S$ is a smooth projective surface,
then a geometric quotient 
$S/R$ can be nonprojective (but proper)
and if $X$ is a smooth proper 3-fold, 
$X/R$ can be an algebraic space which is not a scheme.

\begin{exmp} \label{nonqpglue.exmp}
1. Let $C,D$ be  smooth projective curves
and $S$ the blow up of $C\times D$ at a point $(c,d)$.
Let $C_1\subset S$ be the birational transform of
$C\times \{d\}$, $C_2:=C\times \{d'\}$ for some $d'\neq d$
and $\p^1\cong E\subset S$ the exceptional curve.

Fix an isomorphism  $\sigma:C_1\cong C_2$.
This generates an equivalence relation $R$ which is the identity
on $S\setminus(C_1\cup C_2)$.
As we see in (\ref{quot.by.R.lowdim}), $S/R$ is a surface of finite type.
Note however that the image of $E$ in $S/R$ is
numerically equivalent to $0$, thus $S/R$ is not quasi projective.
Indeed, let $M$ be any line bundle on $S/R$.
Then $\pi^*M$ is a line bundle on $S$
such that $(C_1\cdot \pi^*M)=(C_2\cdot \pi^*M)$.
Since $C_2$ is numerically equivalent to
$C_1+E$, this implies that $(E\cdot \pi^*M)=0$.

2. Take  $S\cong \p^2$ and $Z:=(x(y^2-xz)=0)$. Fix an isomorphism 
of the line $(x=0)$ and
 the conic $(y^2-xz=0)$ which is the identity on
their intersection. As before, this generates an equivalence relation $R$
which is the identity on their complement.
By (\ref{quot.by.R.lowdim}),
 $\p^2/R$ exists as a scheme but it is not projective.

Indeed, if $M$ is a line bundle on $\p^2/R$
then $\pi^*M$ is a line bundle on $\p^2$ whose degree
on a line is the same as its degree on a conic.
Thus $\pi^*H\cong \o_{\p^2}$  and so $H$ is not ample.

3. Let  
 $S=S_1\amalg S_2\cong \p^2\times \{1,2\}$ be 2 copies of
$\p^2$. Let $E\subset \p^2$ be a smooth cubic.
 For a point  $p\in E$, 
let  $\sigma_p: E\times \{1\}\to E\times \{2\}$
be  the identity composed with translation by $p\in E$.
 As before, this generates an equivalence relation $R$
which is the identity on their complement.

Let $M$ be a line bundle on $S/R$.
Then $\pi^*M|_{S_i}\cong \o_{\p^2}(m_i)$ for some $m_i>0$, and we conclude that
$$
\o_{\p^2}(m_1)|{E}\cong \tau_p^*\bigl(\o_{\p^2}(m_2)|{E}\bigr).
$$
This holds iff $m_1=m_2$ and $p\in E$ is a $3m_1$-torsion point.
Thus the projectivity of $S/R$ depends very subtly on
the gluing map  $\sigma_p$.
\end{exmp}

\begin{exmp} Hironaka's example in 
\cite[App.B.3.4.1]{hartsh} gives a
 smooth, proper variety $X$ 
and  two curves $\p^1\cong C_1\cong  C_2\subset X$ such that
$C_1+C_2$ is homologous to $0$. Let $g:C_1\cong C_2$ be an isomorphism
and $R$ the corresponding equivalence relation.

We claim that there is a no quasi projective open subset
$U\subset X$ which intersects both $C_1$ and $C_2$. 
Assume to the contrary that $U$ is such. Then there is an ample divisor
$H_U\subset U$ which intersects both curves but does not contain either.
Its closure $H\subset X$ is a Cartier divisor
which intersects both curves but does not contain either.
Thus $H\cdot (C_1+C_2)>0$, a contradiction.

This shows that if $p\in X/R$ is on the image of $C_i$
then $p$ does not have any affine open neighborhood
since the preimage of an affine set by a finite morphism is again affine.
Thus $X/R$ is not a scheme.
\end{exmp}

\begin{exmp}\cite{lip-lip}\label{sch.quot.exmp}
Fix a field $k$ and let $a_1,\dots, a_n\in k$
be different elements. Set $A:=k[x,y]/\bigl( \prod_i (x-a_iy)\bigr)$. Then
$Y:=\spec R$
is $n$ lines through the origin. Let $f:X\to Y$ its normalization.
Thus $X=\amalg_i \spec  k[x,y]/(x-a_iy)$.
Note that
$$
 k[x,y]/(x-a_iy)\otimes_A k[x,y]/(x-a_jy)=
\left\{
\begin{array}{l}
 k[x,y]/(x-a_iy)\qtq{if $a_i=a_j$, and}\\
 k \qquad \qquad\qquad\ \  \qtq{if $a_i\neq a_j$.}
\end{array}
\right.
$$
Thus $X\times_YX$ is reduced. It is the union  of the diagonal
$\Delta_X$ and of $f^{-1}(0,0)\times f^{-1}(0,0)$.
Thus $X/\bigl(X\times_YX\bigr)$ is 
 a seminormal scheme which  is isomorphic to the $n$ coordinate axes in $\a^n$.
For $n\geq 3$, it is not isomorphic to $Y$.

One can also get similar examples where $Y$ is integral.
Indeed, let $Y\subset \a^2$ be any plane curve
whose only singularities are ordinary multiple points
and let  $f:X\to Y$ be its normalization.
By the above computations, $X\times_YX$ is reduced
and $X/\bigl(X\times_YX\bigr)$ is 
 the seminormalization of $Y$.

If $Y$ is a reduced scheme with normalization $\bar Y\to Y$.
Then, as we see in (\ref{quot.X/S.finite.lem}),
the geometric quotient $\bar Y/\bigl(\bar Y\times_Y\bar Y\bigr)$
exists. It coincides with the strict closure considered
in \cite{lip-arf}. The curve case was introduced earlier
in \cite{arf}.

The related Lipschitz closure is studied in
\cite{pham} and \cite{lip-lip}.
\end{exmp}

\section{Basic results}\label{basic.res.sec}

In this section we prove some basic existence results
for geometric quotients.

\begin{lem}\label{quot.X/S.finite.lem}
Let $S$ be a Noetherian scheme.
Assume that $X$ is finite over $S$ and let $p_1,p_2:R\rightrightarrows X$
be a finite, set theoretic equivalence relation over $S$.
Then the geometric  quotient  $X/R$ exists.
\end{lem}

Proof.  Since  $X\to S$ is affine,  the categorical quotient is
 the spectrum of the
 $\o_S$-algebra
$$
\ker\Bigl[\o_X\stackrel{p_1^*-p_2^*}{\longrightarrow} \o_R\Bigr].
$$
This kernel is a submodule of the finite $\o_S$-algebra $\o_X$,
hence itself a finite $\o_S$-algebra.
The only question is about the geometric fibers of $X\to (X/R)^{cat}$.
Pick any $s\in S$. Taking kernel commutes with
flat base scheme extensions.
 Thus we may assume
that $S$ is complete, local  with closed point $s$ and 
algebraically closed residue field $k(s)$.
We  need to show that the reduced fiber of
$(X/R)^{cat}\to S$ over $s$ is naturally
isomorphic to $\red X_s/\red R_s$.

If $U\to S$ is any finite map then $\o_{\red U_s}$ is  a 
sum of $m(U)$  copies of $k(s)$ for some $m(U)<\infty$. 
$U$ has $m(U)$ connected components $\{U_i:i=1,\dots, m(U)\}$
and each $U_i\to S$ is finite.
 Thus $U\to S$ uniquely factors as
$$
U\stackrel{g}{\to} \amalg_{m(U)}S\to S
\qtq{such that}
g_s:\red U_s\stackrel{\cong}{\longrightarrow} \amalg_{m(U)}\spec k(s)
$$
is an isomorphism, where $\amalg_mS$ denotes the 
disjoint union of $m$ copies of $S$.

Furthermore, if $\{R_j:j=1,\dots m(R)\}$ are the connected components
of $R$, then $\sigma_1$ maps $R_j$ to some $U_i$ and
we have a commutative diagram
$$
\begin{array}{ccc}
R_j&
\stackrel{\sigma_1}{\longrightarrow}
& U_i\\
\downarrow && \downarrow\\
S& =& S.
\end{array}
$$
Applying this to $X\to S$ and $R\to S$, we obtain
a commutative diagram
$$
\begin{array}{lcl}
\o_S^{m(X)}&
\stackrel{p_1^*(s)-p_2^*(s)}{\longrightarrow}
&\o_S^{m(R)}\\
\ \downarrow && \ \downarrow\\
\o_X&\stackrel{p_1^*-p_2^*}{\longrightarrow}& \o_R.
\end{array}
$$
Thus, for $\tau:={p_1^*(s)-p_2^*(s)}$,  we get a morphism
$$
\ker\Bigl[\o_S^{m(X)}
\stackrel{\tau}{\longrightarrow}
 \o_S^{m(R)}\Bigr]\to  \o_{(X/R)^{cat}}.
$$
The kernel on the left is  $m:=|X_s/R_s|$ copies of $\o_S$, hence
we obtain a factorization
$$
(X/R)^{cat}\to \amalg_mS\to S
\qtq{such that}
\red (X/R)^{cat}_s\to \amalg_m\spec k(s)
$$
is an isomorphism.\qed
\medskip

For later reference, we record the following
straightforward consequence.

\begin{cor}\label{quot.of.sub} 
 Let
$R\rightrightarrows X$
 be a finite, set theoretic equivalence relation
such that $X/R$ exists. Let $Z\subset X$ be a closed
$R$-invariant subscheme.
Then $Z/R|_Z$ exists and
 $Z/R|_Z\to X/R$ is  a finite and 
universal homeomorphism
(\ref{univ.homeo.defn}) onto its image.\qed
\end{cor}

\begin{exmp}\label{subquot.neq.sub}
 Even in nice situations, $Z/R|_Z\to X/R$
need not be a closed embedding, as shown by the following examples.
\smallskip

(\ref{subquot.neq.sub}.1)  Set $X:=\a^2_{xy}\amalg \a^2_{uv}$
and let $R$ be the equivalence relation that
identifies the $x$-axis with the $u$-axis.

Let $Z=(y=x^2)\amalg (v=u^2)$. In $Z/R|_Z$ the two
components intersect at a  node,
  but the image of $Z$ in $X/R$ has  a tacnode.

In this example the problem is clearly caused by
ignoring the scheme structure of $R|_Z$.
As the next example shows, similar phenomena happen
even if $R|_Z$ is reduced.

(\ref{subquot.neq.sub}.2)  Set $Y:=(xyz=0)\subset \a^3$. Let
 $X$ be the normalization of $Y$ and $R:=X\times_YX$.
Set $W:=(x+y+z=0)\subset  Y$
 and let $Z\subset X$  be the preimage of $W$.
As computed in
(\ref{sch.quot.exmp}),  $R$ and $R|_Z$ are both reduced, $Z/R|_Z$ is the
seminormalization of $W$ and $Z/R|_Z\to W$ is not an isomorphism.
\end{exmp}

\begin{rem}\label{bb.count.exmp}
 The following counter example to (\ref{quot.X/S.finite.lem})
is proposed in \cite[6.2]{bial}. Consider the  diagram
$$
\begin{array}{ccc}
 \spec k[x,y] & \stackrel{p_1}{\longrightarrow} & \spec k[x,y^2,y^3]\\
p_2\downarrow\hphantom{p_2} && \hphantom{q_2}\downarrow q_2\\
\spec k[x+y,x+x^2, y^2,y^3]&  \stackrel{q_1}{\longrightarrow} &
\spec k[x+x^2, xy^2, xy^3,  y^2,y^3]
\end{array}
\eqno{(\ref{bb.count.exmp}.1)}
$$
It is easy to see that the $p_i$ are homeomorphisms
but $q_2p_1=q_1p_2$ maps  $(0,0)$ and $(-1,0)$ to the same point.
If (\ref{bb.count.exmp}.1)
were a universal push-out, one would get a 
counter example to (\ref{quot.X/S.finite.lem}). 
However, it
is not a universal push-out. Indeed, 
$$
\begin{array}{rcl}
\tfrac13(x+y)^3+\tfrac12(x+y)^2&=&\ \ 
\bigl(\tfrac13 x^3+\tfrac12x^2\bigr)+ (x^2+x)y\ \ \ \ \ \ \ \ 
+xy^2+\tfrac12y^2+\tfrac13 y^3\\[1ex]
&=&
-\bigl(\tfrac23 x^3+\tfrac12x^2\bigr)+ (x^2+x)(x+y)+xy^2+\tfrac12
y^2+\tfrac13 y^3
\end{array}
$$
shows that $\tfrac23 x^3+\tfrac12x^2$ is also in the intersection 
$k[x,y^2,y^3]\cap k[x+y,x+x^2, y^2,y^3]$.
\end{rem}

\medskip
Another case where 
 $X/R$ is easy to obtain is the following.

\begin{lem}\label{quot.pure.dim.lem}
Let $p_1,p_2:R\rightrightarrows X$ be
 a finite, set theoretic equivalence relation where
 $X$ is normal, Noetherian  and $X,R$ are both pure dimensional.
Assume  that one of the following holds:
\begin{enumerate}
\item  $X$ is  defined over a field of characteristic 0,
\item  $X$ is essentially of finite type over $S$, or
\item  $X$  is  defined over a field of characteristic $p>0$ and
 the Frobenius map
$F^p:X\to X^{(p)}$ (\ref{frob.say}) is  finite.
\end{enumerate}
Then the geometric quotient $X/R$ exists as an algebraic space.
$X/R$ is  normal, Noetherian  and essentially of finite type over $S$
in  case (2).
\end{lem}

Proof. Let us first deal with the case when $R$
comes from a  finite group acting on $X$.
This case is a result of Deligne, discussed in \cite[IV.1.8]{Knutson71}.

For  $x\in X$, let $G_x\subset G$ denote the stabilizer.
Let $x\in U_x\subset X$ be a $G_x$-invariant affine open
subset. By shrinking $U_x$
we may assume that  $G_{x'}\subset G_x$ for every $x'\in U_x$.

In the affine case, quotients by finite
groups are easy to get (\ref{inv.of.fin.gps});
this is where the conditions (1--3) are used.
Thus $U_x/G_x$ exists and it is easy to see that the
 $U_x/G_x$ give \'etale charts for $X/G$.

 In  the general  case, it is enough to construct the quotient  when  
 $X$ is irreducible.
Let $m$ be the separable degree of the projections $\sigma_i:R\to X$.

Consider the $m$-fold product
$X\times\cdots\times X$ with coordinate projections $\pi_i$.
Let $R_{ij}$ (resp.\ $\Delta_{ij}$)
denote the preimage of $R$ (resp.\ of the diagonal)
under  $(\pi_i,\pi_j)$. A geometric point of  $\cap_{ij} R_{ij}$
is a sequence of geometric  points
$(x_1,\dots, x_m)$ such that any 2 are $R$-equivalent
and a geometric point of  $\cap_{ij} R_{ij}\setminus \cup_{ij}\Delta_{ij}$
is a sequence $(x_1,\dots, x_m)$ that constitutes a whole
$R$-equivalence class. 
Let $X'$ be the
normalization of the closure of 
$\cap_{ij} R_{ij}\setminus \cup_{ij}\Delta_{ij}$.
Note that every $\pi_{\ell}:\cap_{ij} R_{ij}\to X$ is finite,
hence the projections
$\pi'_{\ell}:X'\to X$ are finite.

The symmetric group $S_m$ acts on $X\times\cdots\times X$
by permuting the factors and this lifts to an
$S_m$-action on $X'$.  Over a dense open subset of $X$,
the $S_m$-orbits on  the  geometric points of $X'$
are exactly the $R$-equivalence classes.

 Let $X^*\subset X'/S_m\times X$ be the image of $X'$ under the
diagonal map.

By construction, 
$X^*\to X$ is finite and one-to-one on geometric points
over an open set. Since $X$ is normal, $X^*\cong X$
in characteristic 0 and $X^*\to X$ is purely inseparable
 in positive characteristic.

In characteristic 0, we thus have a morphism
$X\to X'/S_m$ whose geometric fibers are
exactly the $R$-equivalence classes.
Thus $X'/S_m=X/R$.

Essentially the same works in positive characteristic,
see Section \ref{pos.char.sec} for details. \qed

\begin{lem}\label{quot.norm.lem}
Let $p_1,p_2:R\rightrightarrows X$ be
 a finite, set theoretic equivalence relation
such that $(X/R)^{cat}$ exists.
\begin{enumerate}
\item 
If $X$ is normal and $X, R$ are  pure dimensional 
then  $(X/R)^{cat}$ is also normal.
\item If
$X$ is seminormal
then  $(X/R)^{cat}$ is also seminormal.
\end{enumerate} 
\end{lem}

Proof. In the first case,
let $Z\to (X/R)^{cat}$ be a finite morphism which is an isomorphism
at all generic points of $(X/R)^{cat}$. Since $X$ is normal, 
$\pi:X\to (X/R)^{cat}$ lifts to $\pi_Z:X\to Z$.
 By assumption, 
$\pi_Z\circ p_1$ equals $\pi_Z\circ p_2$ at
all generic points of $R$ and $R$ is reduced. Thus
$\pi_Z\circ p_1=\pi_Z\circ p_2$. 
The universal property of categorical quotients gives
$(X/R)^{cat}\to Z$, hence $Z=(X/R)^{cat}$ and  $(X/R)^{cat}$ is normal. 

In the second case,
let $Z\to (X/R)^{cat}$ be a finite morphism which is a
universal homeomorphism (\ref{univ.homeo.defn}). As before, 
we get liftings
$\pi_Z\circ p_1,\pi_Z\circ p_2:R\rightrightarrows X\to Z$
which agree on closed points. Since $R$ is reduced,
we conclude that $\pi_Z\circ p_1=\pi_Z\circ p_2$, thus 
$(X/R)^{cat}$ is  seminormal.
\qed
\medskip

The following result goes back at least to E.\ Noether.

\begin{prop}\label{inv.of.fin.gps}
 Let $A$ be a Noetherian ring,
$R$ a  Noetherian $A$-algebra and $G$ a finite group of
$A$-automorphisms of $R$. Let $R^G\subset R$ denote the subalgebra of
$G$-invariant elements. Assume  that one of the following holds:
\begin{enumerate}
\item   $\frac1{|G|}\in A$,
\item   $R$ is essentially of finite type over $A$, or
\item  $R$ is finite over
$A[R^p]$ for every prime $p$ that divides $|G|$.
\end{enumerate}
 Then
 $R^G$ is Noetherian and $R$ is finite over $R^G$.
\end{prop}

Proof. Assume first that $R$ is a localization of a
finitely generated $A$ algebra  $A[r_1,\dots, r_m]\subset R$.
We may assume that $G$ permutes the $r_j$.
Let $\sigma_{ij}$ denote the $j$th elementary symmetric polynomial
of the $\{g(r_i):g\in G\}$. Then
$$
A\bigl[\sigma_{ij}\bigr]\subset A[r_1,\dots, r_m]^G\subset R^G
$$
and, with $n:=|G|$,  each $r_i$ satisfies the equation
$$
r_i^n-\sigma_{i1}r_i^{n-1}+\sigma_{i2}r_i^{n-2}-+\cdots =0.
$$
Thus $A[r_1,\dots, r_m]$ is integral over $A\bigl[\sigma_{ij}\bigr]$,
hence also over the larger ring 
$A[r_1,\dots, r_m]^G$.

By assumption $R=U^{-1}A[r_1,\dots, r_m]$
where $U$ is a subgroup of units in $A[r_1,\dots, r_m]$.
We may assume that $U$ is $G$-invariant.
If $r/u\in R$ where $r\in A[r_1,\dots, r_m]$ and
$u$ a unit in $A[r_1,\dots, r_m]$, then
$$
\frac{r}{u}=\frac{r\prod_{g\neq 1} g(u)}{u\prod_{g\neq 1} g(u)},
$$
where the product is over the non-identity elements of $G$.
Thus $r/u=r'/u'$ where $r'\in A[r_1,\dots, r_m]$ and
$u'$ is a $G$-invariant  unit in $A[r_1,\dots, r_m]$.
Therefore, 
$$
R=\bigl(U^G\bigr)^{-1}A[r_1,\dots, r_m]
\qtq{is finite over}
\bigl(U^G\bigr)^{-1}A\bigl[\sigma_{ij}\bigr].
$$
Since $R^G$ is an $\bigl(U^G\bigr)^{-1}A\bigl[\sigma_{ij}\bigr]$-submodule of
$R$, it is also finite  over $\bigl(U^G\bigr)^{-1}A\bigl[\sigma_{ij}\bigr]$,
hence the localization of a  finitely generated algebra.

Assume next that $|G|$ is invertible in $A$.
We claim that $JR\cap R^G=J$ for any ideal $J\subset R^G$.
Indeed, if $a_i\in R^G$, $r_i\in R$ and
$\sum r_ia_i\in R^G$ then
$$
|G|\cdot \sum_i r_ia_i=\sum_{g\in G}\sum_ig(r_i)g(a_i)=
\sum_ia_i\sum_{g\in G}g(r_i)\in \sum_ia_iR^G.
$$
If $|G|$ is invertible,  this gives that
$R^G\cap\sum a_iR =\sum a_iR^G$.
Thus the map $J\mapsto JR$ 
from the ideals of $R^G$ to the ideals of $R$ is an injection which
preserves inclusions. Therefore
 $R^G$ is Noetherian if $R$ is.

If $R$ is an integral domain, then $R$ is  finite over $R^G$
by (\ref{r/rg.finite.general}). The general case,
which we do not use,  is left to the reader.

The arguments in case (3) are quite involved, see \cite{fogarty}.
\qed

\begin{lem} \label{r/rg.finite.general}
Let $R$ be an integral domain and $G$ a finite group
of automorphisms of $R$. Then $R$ is contained in a finite
$R^G$-module. Thus, if $R^G$ is Noetherian, then
$R$ is  finite over $R^G$.
\end{lem}

Proof. Let $K\supset R$ and $K^G\supset R^G$ denote the quotient fields.
 $K/K^G$ is a Galois extension with group $G$.
Pick $r_1,\dots, r_n\in R$ that form a $K^G$-basis of $K$.
Then any $r\in R$ can be written as
$$
r=\sum_i a_i r_i\qtq{where $a_i\in K^G$.}
$$
Applying any $g\in G$ to it, we get a system of equations
$$
\sum_i g(r_i)a_i= g(r)\qtq{for $g\in G$.}
$$
We can view these as  linear equations with unknowns $a_i$.
The system determinant is $D:=\det_{i,g}\bigl(g(r_i)\bigr)$,
which is nonzero since its square is the discriminant of
 $K/K^G$. $D$ is $G$-invariant up to sign, thus
$D^2$ is  $G$-invariant hence in $R^G$.
By Kramer's rule, $a_i\in D^{-2}R^G$, hence
$R\subset D^{-2}\sum_i r_i R^G$.\qed

\medskip
In the opposite case, when the equivalence relation 
is nontrivial only on a proper subscheme, we have the
following general result.

\begin{prop}\label{quot.by.pushout} Let $X$ be  a reduced scheme,
$Z\subset X$ a closed, reduced subscheme and
 $R\rightrightarrows X$
 a finite, set theoretic equivalence relation.
Assume that $R$ is the identity on $R\setminus Z$
and that the geometric quotient  $Z/R|_Z$ exists.
Then $X/R$ exists and is given by the
universal push-out diagram
$$
\begin{array}{ccc}
Z & \into & X\\
 \downarrow && \downarrow \\
Z/R|_Z & \into & X/R.
\end{array}
$$
\end{prop}

Proof. Let $Y$ denote the universal push-out (\ref{glue.thm.asp}).
Then $X\to Y$ is finite and so $X/R$ exists 
and we have a natural map $X/R\to Y$ by
(\ref{quot.X/S.finite.lem}).
 On the other hand, there is a natural map $Z/R|_Z\to X/R$
by (\ref{quot.of.sub}), hence the universal property of the push-out gives
the inverse $Y\to X/R$.\qed

\section{Inductive plan for constructing quotients}\label{induct.plan.sect}

\begin{defn} \label{up-down.eq.defn}
Let 
$R\rightrightarrows X$
be a finite, set theoretic equivalence relation
and $g:Y\to X$ a finite morphism. Then
$$
g^*R:= R\times_{(X\times X)}(Y\times Y)\rightrightarrows Y
$$
defines a finite, set theoretic equivalence relation on $Y$.
It is called the {\it pull-back} of $R\rightrightarrows X$.
(Strictly speaking, it should be denoted by
$(g\times g)^*R$.)

Note that the $g^*R$-equivalence classes on the 
geometric points of $Y$  map injectively to
 the $R$-equivalence classes on the 
geometric points of $X$.

If $X/R$  exists then, by
(\ref{quot.X/S.finite.lem}),  $Y/g^*R$ also exists and the
natural morphism $Y/g^*R\to X/R$  is
injective on geometric points. If, in addition, $g$ is
surjective then $Y/g^*R\to X/R$ is a finite and 
universal homeomorphism
(\ref{univ.homeo.defn}). Thus, if $X$ is seminormal
and the characteristic is 0, then
$Y/g^*R\cong X/R$.

Let $h:X\to Z$ be a finite morphism.
If the geometric fibers of $h$
are subsets of $R$-equivalence classes, then the composite
$R\rightrightarrows X\to Z$
defines a finite, set theoretic pre-equivalence relation 
$$
h_*R:=(h\times h)(R)\subset Z\times Z,
$$
called the {\it push forward} of $R\rightrightarrows X$.
If  $Z/R$  exists, then, by (\ref{quot.X/S.finite.lem}),
 $X/R$ also exists and the
natural morphism  $X/R\to Z/R$ is a finite and 
universal homeomorphism.
\end{defn}

\begin{lem}\label{norm.quot=quot}
 Let $X$ be weakly normal, excellent and
$R\rightrightarrows X$
 a finite, set theoretic equivalence relation.
Let $\pi:X^n\to X$ be the normalization and
$R^n\rightrightarrows X^n$ the pull back of $R$ to $X^n$.
If $X^n/R^n$ exists then $X/R$ also exists
and $X/R=X^n/R^n$.
\end{lem}

Proof. 
Let $X^*\subset \bigl(X^n/R^n\bigr)\times_S X$ be the image of
$X^n$ under the diagonal morphism.
 Since $X^n\to X$ is a finite
surjection, $X^*$ is a closed subscheme of
 $\bigl(X^n/R^n\bigr)\times_S X$ and
$X^*\to X$ is  a finite
surjection. Moreover, for any geometric point $\bar x\to X$,
its preimages $\bar x_i\to X^n$ are $R^n$-equivalent, hence
they map to the same point in $\bigl(X^n/R^n\bigr)\times_S X$.
Thus $X^*\to X$ is  a  finite and 
one-to-one on geometric points, so it is a  finite and 
universal homeomorphism
(\ref{univ.homeo.defn}). $X^n\to X$ is a local isomorphism at
 the generic point of every irreducible component of $X$,
hence $X^*\to X$ is also a local isomorphism at
 the generic point of every irreducible component of $X$.
Since $X$ is weakly normal,  $X^*\cong X$ and
we have a morphism $X\to X^n/R^n$ and thus  $X/R=X^n/R^n$.\qed

\begin{lem}\label{normal.subrel.lem}
Let $X$  be normal and of pure dimension $d$.
  Let
$\sigma:R\rightrightarrows X$ be
 a finite, set theoretic equivalence relation
and $R^d\subset R$  its $d$-dimensional part.
Then $\sigma^d:R^d\rightrightarrows X$ is also an equivalence relation.
\end{lem}

Proof. The only question is transitivity.
Since $X$ is normal, the maps $\sigma^d:R^d\rightrightarrows X$
are both universally open  by Chevalley's criterion, see
\cite[IV.14.4.4]{EGA}.
Thus the fiber product $R^d\times_XR^d\to X$ is also
universally open and hence its irreducible
components have  pure dimension $d$.\qed

\begin{exmp} Let $C$ be a curve with an involution
$\tau$. Pick $p,q\in C$ with $q$ different from $p$  and $\tau(p)$.
Let $C'$ be the nodal curve obtained from $C$ by identifying $p$ and $q$.
The equivalence relation generated by $\tau$ on $C'$
consists of the diagonal, the graph of $\tau$ plus the pairs
$\bigl(\tau(p),\tau(q)\bigr)$ and $\bigl(\tau(q),\tau(p)\bigr)$.
The 1-dimensional parts of the equivalence relation
do not form an equivalence relation.
\end{exmp}

\begin{say}[Inductive plan] \label{induct.plan}
Let $X$ be an excellent scheme 
that satisfies one of the conditions
(\ref{quot.pure.dim.lem}.1--3)
and 
$R\rightrightarrows X$
 a finite, set theoretic equivalence relation.
We aim to construct the
 geometric quotient $X/R$ in two steps.
First  we construct a space that, roughly speaking,  should be the
normalization of  $X/R$ and then
 we try to go  from the normalization to the geometric quotient itself.
\medskip

{\it Step 1.} Let $X^n\to X$ be the normalization of $X$ and
$R^n\rightrightarrows X^n$ the pull back of $R$ to $X^n$.
Set $d=\dim X$ and let $X^{nd}\subset X^n$ 
(resp.\ $R^{nd}\subset R^n$) denote the union of the 
$d$-dimensional irreducible components. 
By  (\ref{normal.subrel.lem}),  $R^{nd}\rightrightarrows X^{nd}$
is a pure dimensional, finite, set theoretic equivalence relation
and the geometric quotient
 $X^{nd}/R^{nd}$ exists by (\ref{quot.pure.dim.lem}).

There is a closed, reduced  subscheme $Z\subset X^n$ of 
dimension $<d$ such that
$Z$ is closed under $R^n$ and the two equivalence relations
$$
R^n|_{X^n\setminus Z}\qtq{and}
R^{nd}|_{X^n\setminus Z} \qtq{coincide.}
$$
Let $Z_1\subset  X^{nd}/R^{nd}$ denote the image of
$Z$. 
 $R^n|_Z\rightrightarrows Z$ gives a
finite set theoretic equivalence relation on $Z$.
Since the geometric fibers of $Z\to Z_1$
are subsets of $R^n$-equivalence classes, by (\ref{up-down.eq.defn}),
 the composite maps  
$R^n|_Z\rightrightarrows Z\to Z_1$
define a
finite set theoretic pre-equivalence relation on $Z_1$.
\medskip

{\it Step 2.} In order to go from $  X^{nd}/R^{nd}$ to  $X/R$,
we  make the following
\medskip

\noindent {\it Inductive assumption} (\ref{induct.plan}.2.1). 
The geometric quotient   $Z_1/\bigl(R^{n}|_Z\bigr)$ exists.
\medskip

Then, by (\ref{quot.by.pushout}) $X^n/R^n$ exists and is given as 
the  universal push-out of the following diagram:
$$
\begin{array}{ccc}
Z_1 & \into & X^n/R^{nd}\\
\downarrow &&\downarrow \\
Z_1/\bigl(R^{n}|_Z\bigr) &\into & X^n/R^n.
\end{array} 
$$

As in (\ref{norm.quot=quot}),
let $X^*\subset \bigl(X^n/R^n\bigr)\times_S X$ be the image of
$X^n$ under the diagonal morphism.
We have established that 
 $X^*\to X$ is  a  finite and 
universal homeomorphism
(\ref{univ.homeo.defn})
sitting in 
 the following diagram:
$$
\begin{array}{ccccccc}
Z_1 & \into & X^n/R^{nd} & \leftarrow & X^n &&\\
\downarrow && \downarrow &\swarrow & \downarrow & \searrow &\\
Z_1/\bigl(R^{n}|_Z\bigr) & \to & X^n/R^n & \leftarrow & X^* & \to & X
\end{array}
\eqno{(\ref{induct.plan}.2.2)}
$$

There are now two ways to proceed.
\medskip

{\it Positive characteristic} (\ref{induct.plan}.2.3).   
Most finite and universal homeomorphisms can be
inverted, up to a power of the Frobenius (\ref{factor.through.frob.prop}),
 and so we obtain a morphism
$$
X\to \bigl(X^*\bigr)^{(q)}\to \bigl(X^n/R^n\bigr)^{(q)}
$$
 for some $q=p^m$.
$X/R$ is then obtained using (\ref{quot.X/S.finite.lem}).
  This is discussed in Section \ref{pos.char.sec}.

In this case the inductive assumption
(\ref{induct.plan}.2.1)  poses no extra problems.
\medskip

{\it Zero characteristic} (\ref{induct.plan}.2.4).   
As the examples of Section \ref{first.exmp.sec} show,
finite and universal homeomorphisms cause a substantial problem.
The easiest way to overcome these difficulties is to assume
to start with that $X$ is seminormal. 
In this case, by (\ref{norm.quot=quot}), we obtain that 
$X/R=X^n/R^n$.

Unfortunately,   the inductive assumption (\ref{induct.plan}.2.1)
 becomes quite restrictive.
By construction $Z_1$ is reduced, but it need not be seminormal
in general. 
Thus we get the induction going only if we can guarantee
that  $Z_1$ is seminormal.
 Note that, because of the inductive
set-up, seminormality needs to hold not only for $X$ and  $Z_1$,
 but on further schemes that one obtains in applying  the
inductive proof to $R^{n}|_Z\rightrightarrows Z_1$, and so on.

It turns out, however,  that  the above inductive plan works
when gluing semi-log-canonical schemes.
This will be treated elsewhere.

\end{say}

\begin{defn}\label{monom.defn} A morphism of schemes
$f:X\to Y$ is a {\it monomorphism} if for every scheme $Z$
the induced map of sets
$\mor(Z,X)\to \mor(Z,Y)$ is an injection.

By \cite[IV.17.2.6]{EGA} this is equivalent to assuming that
$f$ be universally injective and unramified.

A proper monomorphism $f:Y\to X$ is a closed embedding.
Indeed, a proper monomorphism is injective
on geometric points, hence finite. 
Thus it is a closed embedding
iff $\o_X\to f_*\o_Y$ is onto. By the Nakayama lemma this is
equivalent to
$f_x:f^{-1}(x)\to x$ being an isomorphism for every
 $x\in f(Y)$. By passing to geometric points,
we are down to the case when 
 $X=\spec k$, $k$ is algebraically closed and
  $Y=\spec A$ where $A$ is an Artin $k$-algebra. 

If $A\neq k$ then there are at least 2 different $k$ maps
$A\to k[\epsilon]$, thus $\spec A\to \spec k$ is
not a monomorphism.
\end{defn}

\begin{defn}\label{univ.homeo.defn}
 We say that a morphism of schemes $g:U\to V$ is a {\it
universal homeomorphism} if  it is a homeomorphism and for every $W\to V$ the
induced morphism $U\times_VW\to W$ is again a homeomorphism.
The definition extends to morphisms of algebraic spaces the usual way
\cite[II.3]{Knutson71}.

A simple example of a  homeomorphism which is not a 
universal homeomorphism is $\spec K\to \spec L$ where $L/K$ is a finite
field extension and $L\neq K$. A more interesting example is
given by the normalization of
the nodal curve $\bigl(y^2=x^2(x+1)\bigr)$ with one of 
the preimages of the node
removed:
$$
\a^1\setminus \{-1\}\to \bigl(y^2=x^2(x+1)\bigr)
\qtq{given by} t\mapsto \bigl(t^2-1, t(t^2-1)\bigr).
$$
When  $g$ is finite, the notion is pretty much
set theoretic
since a continuous proper map of topological spaces
which is injective and surjective
is a homeomorphism. Thus we see that for a finite and surjective morphism of
algebraic spaces
 $g:U\to V$ the following are equivalent
(cf.\  \cite[I.3.7--8]{ega71})
\begin{enumerate}
\item $g$ is  a universal homeomorphism.
\item $g$ is surjective and universally injective.
\item For every $v\in V$ the fiber  $g^{-1}(v)$ has a single 
point $v'$ and 
$k(v')$ is a purely inseparable field extension of   $k(v)$.
\item $g$ is surjective and injective on geometric points.
\end{enumerate}

One of the most important properties of these morphisms is that
taking the fiber product induces an
equivalence between the categories 
$$
(\mbox{\'etale morphisms: $\ast\to V$})
\stackrel{\ast\mapsto \ast\times_VU}{\longrightarrow}
(\mbox{\'etale morphisms: $\ast\to U$}).
$$
See \cite[IX.4.10]{sga1} for a proof. 
We do not use this in the sequel.
\end{defn}

In low dimensions one can start the  method (\ref{induct.plan}) and
it gives the following. These results are sufficient to deal with the
moduli problem for surfaces.

\begin{prop}\label{quot.by.R.lowdim} Let $S$ be a Noetherian scheme 
over a field of characteristic 0 and
$X$ an algebraic space of finite type over $S$.
Let $R\rightrightarrows X$ be a finite, set theoretic equivalence relation.
Assume that one of the following holds: 
\begin{enumerate}
\item $X$ is  1-dimensional and reduced, 
\item $X$ is  2-dimensional and seminormal,
\item $X$ is  3-dimensional, normal and there is a closed, seminormal
$Z\subset X$ such that $R$ is the identity on $X\setminus Z$. 
\end{enumerate}
Then the geometric  quotient $X/R$ exists.
\end{prop}

Proof.
Consider first the case when $\dim X=1$.
Let $\pi:X^n\to X$ be the normalization.
We construct  $X^n/R^{nd}$ as in (\ref{induct.plan}).
Note that since $Z$ is zero dimensional, it is finite over $S$.
Let $V\subset S$ be its image.
Next we make a different choice for $Z_1$.
Instead, we take a subscheme $Z_2\subset X^n/R^{nd}$
whose support is $Z_1$
such that the pull back of its ideal sheaf $I(Z_2)$ 
to $X^n$ is a subsheaf of the inverse image sheaf  
$\pi^{-1}\o_X\subset \o_{X^n}$.

Then we consider the  push-out  diagram
$$
V\leftarrow Z_2\into X^n/R^{nd}
$$
 with universal push-out  $Y$.
Then $X\to Y$ is a finite  morphism and 
 $X/R$ exists by (\ref{quot.X/S.finite.lem}).

The case when $\dim X=2$ and $X$ is seminormal
is a direct consequence of 
(\ref{induct.plan}.2.4) since the inductive assumption
(\ref{induct.plan}.2.1) is guaranteed by (\ref{quot.by.R.lowdim}.1).

If $\dim X=3$, then $X$ is already normal and
 $Z$ is seminormal by assumption. Thus
$Z/\bigl(R|_Z\bigr)$ exists by (\ref{quot.by.R.lowdim}.2).
Therefore $X/R$ is given by the push-out of
$Z/\bigl(R|_Z\bigr)\leftarrow Z\into X$.
\qed

\section{Quotients in positive characteristic}\label{pos.char.sec}

The main result of this section is the proof of
(\ref{quot.by.R.charp}).

\begin{say}[Geometric Frobenius Morphism] \cite[XIV]{sga5}
\label{frob.say}
Let $S$ be an $\f_p$-scheme.
Fix $q=p^r$ for some natural number $r$. Then $a\mapsto a^q$ defines an
$\f_p$-morphism  $F^q:S\to S$. This can be extended to polynomials
by the formula
$$
f=\sum a_Ix^I\mapsto f^{(q)}:= \sum a^q_Ix^I.
$$

Let $U=\spec R$ be an affine scheme over $S$. Write
$R=\o_S[x_1,\dots,x_m]/(f_1,\dots,f_n)$ and set
$$
R^{(q)}:= \o_S[x^{(q)}_1,\dots,x^{(q)}_m]/(f^{(q)}_1,\dots,f^{(q)}_n)
\qtq{and}  U^{(q)}:=
\spec R^{(q)},
$$
where the $x^{(q)}_i$ are new variables.
There are natural morphisms
$$
F^q:U\to U^{(q)}\qtq{and} (F^q)^*:R^{(q)}\to R\qtq{given by}
(F^q)^*(x^{(q)}_i)=x_i^q.  
$$

It is easy to see that these are independent of the choices made. Thus $F^q$
gives a functor from algebraic spaces over $S$  to algebraic spaces over
 $S$.
One can define $X^{(q)}$ intrinsically as 
$$
X^{(q)}=X\times_{S,F^q} S.
$$

 If  $X$ is an algebraic space 
which is essentially of finite type over $\f_p$
then
$F^q:X\to X^{(q)}$ is a finite and universal homeomorphism.
\end{say}

For us the most important feature of the Frobenius morphism is the
following universal property:

\begin{prop}\label{factor.through.frob.prop}
Let $S$ be a scheme essentially of finite type over $\f_p$ and 
  $X,Y$    algebraic spaces  which are essentially of finite type over  $S$.
Let $g:X\to Y$ be a finite  and universal homeomorphism. Then for $q=p^r\gg 1$
the map   $F^q$ can be factored as
$$
F^q:X\stackrel{g}{\to}Y \stackrel{\bar g}{\to} X^{(q)}.
$$

Moreover, for large enough $q$ (depending on $g:X\to Y$), 
there is a functorial choice of the 
 factorization  in the sense that if 
$$
\begin{array}{ccc}
X_1& \stackrel{g_1}{\to} & Y_1\\
\downarrow && \downarrow\\
X_2& \stackrel{g_2}{\to} & Y_2\\
\end{array}
$$
is a commutative diagram where the  $g_i$ are 
finite and universal homeomorphisms, then, for $q\gg 1$
(depending on the $g_i:X_i\to Y_i$)  the factorization gives a
commutative diagram
$$
\begin{array}{ccccc}
X_1& \stackrel{g_1}{\to} & Y_1& \stackrel{\bar g_1}{\to} & X_1^{(q)}\\
\downarrow && \downarrow && \downarrow\\
X_2& \stackrel{g_2}{\to} & Y_2& \stackrel{\bar g_2}{\to} & X_2^{(q)}.
\end{array}
$$
\end{prop}

Proof.  It is sufficient
to construct the functorial choice of the 
 factorization in case $X$ and $Y$ are affine schemes over an affine scheme
$\spec C$.
Thus we have a ring homomorphism $g^*:A\to B$ where
$A$ and $B$ are finitely generated $C$-algebras.
We can decompose $g^*$ into $A\twoheadrightarrow B_1$ and $B_1\hookrightarrow
B$. We deal with them separately.

First consider $B_1\subset B$. 
In this case there is no choice involved and
we need to show that there is a $q$ such that
$B^q\subset B_1$, where $B^q$ denotes the $C$-algebra generated by the $q$-th
powers of all elements. The  proof is by Noetherian induction. 

First consider the case when  $B$ is Artinian.
The residue field of $B$ is finite and purely inseparable
over the  residue field of $B_1$, hence $B^q$ is contained in a
field of representatives of $B_1$ for large enough $q$.
 
In the general, we can use the Artinian case over the
generic points to obtain that  $B_1\subset B_1B^q$ is an isomorphism at all
generic points for
$q\gg 1$. Let $I\subset B_1$ denote the conductor of this extension. 
That is, $IB_1B^q=I$.
By induction we know that there is a $q'$ such that 
$(B_1B^q/I)^{q'}\subset B_1/I$. Thus we get that
$$
B^{(qq')}\to B^{qq'}\subset (B_1B^q)^{q'}\subset B_1+IB_1B^q=B_1.
$$

Next consider $A\twoheadrightarrow  B_1$.  
Here we have to make a good choice.
The kernel
 is a nilpotent ideal $I\subset A$, say $I^m=0$. Choose $q'$ such that $q'\geq
m$. For $b_1\in B_1$ let $b_1'\in A$ be any preimage. Then
$(b'_1)^{q'}$ depends only on $b_1$. The map
$$
b_1\mapsto (b'_1)^{q'}\qtq{defines a factorization} 
B_1^{(q')}\to A\to B_1.
$$
Combining the map  $B^{(q)}\to B_1$ with  $B_1^{(q')}\to A$ we
obtain  $B^{(qq')}\to A$.\qed

\begin{say}[Proof of (\ref{quot.by.R.charp})]
The question is local on $S$, hence we may assume that
$S$ is affine. 
$X$ and $R$ are defined over a finitely generated subring
of $\o_S$, hence we may assume that
$S$ is of finite type over $\f_p$.

 The proof is by induction on $\dim X$.
We follow the inductive plan in (\ref{induct.plan})
and use its notation.

If $\dim X=0$ then $X$ is finite over $S$ and 
the assertion follows from (\ref{quot.X/S.finite.lem}).

In going from dimension $d-1$ to $d$, the
 assumption (\ref{induct.plan}.2.1) holds by induction.
Thus (\ref{induct.plan}.2.3) shows that $X^n/R^n$ exists.

Let $X^*\subset (X^n/R^n)\times_S X$ be the image of
$X^n$ under the diagonal morphism. 
As we noted in (\ref{induct.plan}),
 $X^*\to X$ is  a finite and universal homeomorphism.
Thus, by (\ref{factor.through.frob.prop}),
there is a factorization
$$
X^*\to X\to {X^*}^{(q)}\to \bigl(X^n/R^n\bigr)^{(q)}.
$$
Here $X\to \bigl(X^n/R^n\bigr)^{(q)}$ is finite
and $R$ is an equivalence relation on $X$ over the base scheme 
 $\bigl(X^n/R^n\bigr)^{(q)}$.
 Hence, by
(\ref{quot.X/S.finite.lem}), the geometric quotient
$X/R$ exists.\qed
\end{say}

\begin{rem} 
Some of the scheme theoretic aspects of 
the purely inseparable case are treated in
\cite{eke} and 
\cite[Exp.V]{sga3}.
\end{rem}

\section{Gluing or Pinching}\label{glue.sec}

The aim of this section is to give an elementary proof of
the following.

\begin{thm} \cite[Thm.3.1]{artin70}
\label{glue.thm.asp}
Let $X$ be a Noetherian  algebraic space   over
a Noetherian base scheme $S$.
Let $Z\subset X$ be a closed subspace. Let
$g:Z\to V$ be a finite surjection. Then there is a universal
push-out diagram of algebraic spaces
$$
\begin{array}{ccl}
Z & \into & X\\
g \downarrow\hphantom{g} && \ \downarrow \pi\\
V & \into & Y:=X/(Z\to V)
\end{array}
$$
Furthermore,
\begin{enumerate}
\item  $Y$ is a  Noetherian  algebraic space over  $S$
\item  $V\to Y$ is a closed embedding  and $Z=\pi^{-1}(V)$,
\item the natural map 
$\ker\bigl[\o_Y\to \o_V\bigr] \to \pi_*\ker\bigl[\o_X\to \o_Z\bigr]$
is an isomorphism, and
\item if $X$  is of finite type  over $S$ then so is $Y$.
\end{enumerate}
\end{thm}

\begin{rem} If  $X$  is of finite type  over $A$
and  $A$ itself is of finite type  over a field 
or an excellent Dedekind ring, then this is
 an easy consequence of the contraction
results \cite[Thm.3.1]{artin70}. 
The more general case above follows using the
later approximation results \cite{pop}.
The main point of
\cite{artin70} is to understand the case when
$Z\to V$ is proper but not finite. This is much harder than the
finite case we are dealing with.
An elementary approach 
following \cite{ferrand} and \cite{raoult}
is discussed below.
\end{rem}

\begin{say}\label{glue.lem.affine} The affine case of (\ref{glue.thm.asp})
is simple algebra.
Indeed, 
let $R=\o_X$, $I=I(Z)$, $q:\o_X\to \o_Z$ the restriction  and $S=\o_V$.
By (\ref{eak-nag}), $q^{-1}(S)$ is Noetherian.
Set $Y:=\spec q^{-1}(S)$.

 If $\bar r_i\in \o_X/I(Z)$ generate $\o_X/I(Z)$ as an $\o_V$-module
then $r_i\in \o_X$ and $I(Z)$ generate $\o_X$ as a $q^{-1}(\o_V)$-module.
Since $I(Z)\subset q^{-1}(S)$, we obtain that 
 $r_i\in \o_X$ and $1\in \o_X$ generate $\o_X$ as a $q^{-1}(S)$-module.
Applying  (\ref{eak-nag}) to $R_1=\o_X$ and $R_2=q^{-1}(\o_V)$ 
gives the rest.\qed
\end{say}

For the proof of the following result, see
 \cite[Thm.3.7]{mats-cr} and the proof of
(\ref{inv.of.fin.gps}).

\begin{thm}[Eakin-Nagata] \label{eak-nag}
Let  $R_1\supset R_2$ be $A$-algebras
with $A$  Noetherian.
Assume that  $R_1$ is finite over  $R_2$.
\begin{enumerate}
\item If   $R_1$ is Noetherian then so is $R_2$.
\item If 
$R_1$ is a finitely generated $A$-algebra then so is $R_2$. \qed
\end{enumerate}
\end{thm}

Gluing for  algebraic spaces, following \cite{raoult},
is easier than the  quasi projective  case.

\begin{say}[Proof of (\ref{glue.thm.asp})]\label{pf.of.glue.thm.asp}
For every $p\in V$ we construct a commutative diagram
$$
\begin{array}{ccccc}
V_p & \stackrel{g_p}{\leftarrow} & Z_p & \to & X_p\\
\tau_V\downarrow\hphantom{\tau_V} && 
\hphantom{\tau_Z}\downarrow\tau_Z &&\hphantom{\tau_X}\downarrow\tau_X \\
V & \stackrel{g}{\leftarrow} & Z & \to & X
\end{array}
$$
where 
\begin{enumerate}
\item $V_p, Z_p, X_p$ are affine,
\item $g_p$ is finite and $Z_p  \to  X_p$ is a closed embedding,
\item  $V_p$ (resp.\ $Z_p, X_p$) is an  \'etale neighborhood of $p$
(resp.\  $g^{-1}(p)$) and
\item both squares are fiber products.
\end{enumerate}
Affine gluing (\ref{glue.lem.affine}) then gives 
 $Y_p:=X_p/(Z_p\to V_p)$
and (\ref{glue.etloc.lem}) shows that the $Y_p$
are \'etale charts on  $Y=X/(Z\to V)$.

Start with   affine, \'etale neighborhoods $V_1\to V$ of $p$
and $X_1\to X$ of $g^{-1}(p)$. Set $Z_1:=Z\times_XX_1\subset X_1$.
By (\ref{et.nbhds.say})
 we may assume that there is a (necessarily \'etale) morphism
$Z\times_VV_1\to Z_1$. 
In general there is no \'etale neighborhood $X'\to X_1$
extending $Z\times_VV_1\to Z_1$, but there is an affine, \'etale
neighborhood $X_2\to X_1$ extending $Z\times_VV_1\to Z_1$
over a Zariski  neighborhood of $g^{-1}(p)$ (\ref{et.nbhds.say}).

Thus we have  affine, \'etale neighborhoods $V_2\to V$ of $p$,
 $X_2\to X$ of $g^{-1}(p)$ and an open embedding
$ Z\times_XX_2\into Z_2:=Z\times_VV_2$.
Our only remaining problem is that $Z_2\neq Z\times_XX_2$, hence
$Z_2$ is not a subscheme of $X_2$. We achieve this by further shrinking
$V_2$ and $X_2$.

The complement $B_2:=Z_2\setminus Z\times_XX_2$ is closed,
thus $g(B_2)\subset V_2$ is a closed subset not containing $p$.
Pick  $\phi\in \Gamma(\o_{V_2})$ that vanishes on
$g(B_2)$ such that $\phi(p)\neq 0$.
Then $\phi\circ g$ is a function on $Z_2$
that vanishes on $B_2$ but is nowhere zero on $g^{-1}(p)$.
We can thus extend $\phi\circ g$ to a function
$\Phi$ on $X_2$.  Thus
$V_P:=V_2\setminus (\phi=0), Z_P:=Z_2\setminus (\phi\circ g=0)$
and $X_P:=X_2\setminus (\Phi=0)$ have the required properties. \qed
\end{say}

\begin{say}\label{et.nbhds.say} During the proof we have used
two basic properties of \'etale neighborhoods.

First, if
 $\pi:X\to Y$ is finite then for every \'etale neighborhood
$(u\in U)\to (x\in X)$ there is an \'etale neighborhood
$(v\in V)\to (\pi(x)\in Y)$ and a connected
component $(v'\in V')\subset X\times_YV$ such that there is a
lifting  $(v'\in V')\to (u\in U)$.

Second, if $\pi:X\to Y$ is a closed embedding,
$U\to X$ is \'etale and $P\subset U$ is a finite set of points
then we can find an \'etale $V\to Y$
such that $P\subset V$ and there is an open embedding
$(P\subset X\times_YV)\to (P\subset U)$.

For proofs see \cite[3.14 and 4.2--3]{milne}.
\end{say}

The next result shows that gluing commutes with flat morphisms.

\begin{lem} \label{glue.etloc.lem} For $i=1,2$, let
 $X_i$ be  Noetherian affine $A$-schemes,
 $Z_i\subset X_i$  closed subschemes and 
$g_i:Z_i\to V_i$  finite surjections
with universal push-outs $Y_i$.
Assume that in the diagram below  both squares are fiber products.
$$
\begin{array}{ccccc}
V_1 & \stackrel{g_1}{\leftarrow} & Z_1 & \to & X_1\\
\downarrow && \downarrow && \downarrow \\
V_2 & \stackrel{g_2}{\leftarrow} & Z_2 & \to & X_2
\end{array}
$$
\begin{enumerate}
\item If the
vertical maps are flat then $Y_1 \to Y_2$ is also flat.
\item If the
vertical maps are smooth then $Y_1 \to Y_2$ is also smooth.
\end{enumerate}
\end{lem}

Proof. We may assume that all occurring schemes are affine.
Thus we have  $I_i\subset R_i$ and $S_i\subset R_i/I_i$.
Furthermore, $R_1$ is flat over $R_2$,
$I_1=I_2R_1$ and $S_1$ is flat over $S_2$. 
We may also assume that $R_2$ is local.
The key point
is the isomorphism
$$
\bigl( R_1/I_1\bigr)\cong 
\bigl( R_2/I_2\bigr)\otimes_{R_2}R_1\cong \bigl( R_2/I_2\bigr)\otimes_{S_2}S_1.
\eqno{(\ref{glue.etloc.lem}.3)}
$$
Note that this isomorphism is not naturally given,
see (\ref{twist-glue.rem}).

We check the local criterion of flatness
(cf. \cite[Thm.22.3]{mats-cr}). The first condition we need is that
$q_1^{-1}(S_1)/I_1\cong S_1$ be flat over $q_2^{-1}(S_2)/I_2\cong S_2$.
This holds by assumption. Second, we need that the maps
$$
\bigl( I_2^n/I_2^{n+1}\bigr)\otimes_{S_2}S_1\to 
I_2^nR_1/I_2^{n+1}R_1
$$
be isomorphisms.  Since $R_1$ is flat over $R_2$, the
right hand side is isomorphic to
$$
\bigl(I_2^n/I_2^{n+1}\bigr)\otimes_{R_2/I_2} \bigl(R_1/I_1\bigr).
$$
Using (\ref{glue.etloc.lem}.3), we get that
$$
\bigl(I_2^n/I_2^{n+1}\bigr)\otimes_{R_2/I_2} \bigl(R_1/I_1\bigr)\cong
\bigl(I_2^n/I_2^{n+1}\bigr)\otimes_{R_2/I_2} \bigl(R_2/I_2\bigr)
\otimes_{S_2}S_1\cong
\bigl(I_2^n/I_2^{n+1}\bigr)\otimes_{S_2}S_1.
$$
This settles flatness. In order to prove the smooth case,
we just need to check that the fibers of
$Y_1 \to Y_2$ are smooth. Outside $V_1\to V_2$
we have the same fibers as before and $V_1\to V_2$ is
smooth by assumption.\qed

\begin{rem} \label{twist-glue.rem}
Note that there is some subtlety in 
(\ref{glue.etloc.lem}). Consider the simple case
when $X_2$ is a smooth curve over a field $k$,
$Z_2=\{p,q\}$ two $k$-points and $V_2=\spec k$.
Then $Y_2$ is a nodal curve where $p$ and $q$ are identified.

Let now $X_1=X_2\times\{0,1\}$ as 2 disjoint copies. 
Then $Z_1$ consists of 4 points $p_0,q_0,p_1,q_1$
and $V_1$ is 2 copies of $\spec k$.
There are two distinct way to arrange $g_1$.
Namely, 
\begin{enumerate}
\item[--] either $g'_1(p_0)=g'_1(q_0)$ and $g'_1(p_1)=g'_1(q_1)$
and then $Y'_1$ consists of 2 disjoint nodal curves,
\item[--] or $g''_1(p_0)=g''_1(q_1)$ and $g''_1(p_1)=g''_1(q_0)$
and then $Y''_1$ consists of a connected  curve with 2 nodes
and 2 irreducible components.
\end{enumerate}
 Both of these are \'etale double covers
of $Y_2$.
\end{rem}

As in (\ref{pf.of.glue.thm.asp}), 
the next lemma will be used to reduce quasi projective gluing
to the affine case.

\begin{lem} \label{affine.red.lem}
Let $X$ be an $A$-scheme,
 $Z\subset X$ a closed subscheme and
$g:Z\to V$  a finite surjection.

Let $P\subset V$ be a finite subset
and assume that there are open  affine subsets
$P\subset V_1\subset V$ and $g^{-1}(P)\subset X_1\subset X$.

Then there are open  affine subsets
$P\subset V_P\subset V_1$ and $g^{-1}(P)\subset X_P\subset X_1$
such that $g$ restricts to a finite morphism
$g: Z\cap X_P\to V_P$.
\end{lem}

Proof. There is an affine subset
$g^{-1}(P)\subset X_2\subset X_1$ such that
$g^{-1}(V\setminus V_1)$ is disjoint from $X_2$. 
Thus $g$ maps $Z\cap X_2$ to $V_1$.
The problem is that $(Z\cap X_2)\to V_1$ is only quasi finite
in general.
The set $Z\setminus X_2$ is closed in $X$ and so
$g( Z\setminus X_2)$ is closed in $V$.
Since $V_1$ is affine, there is a function
$f_P$ on $V_1$ which vanishes on $g( Z\setminus X_2)\cap V_1$
but does not vanish on $P$. 
Then $f_P\circ g$ is a function on
$g^{-1}(V_1)$ which vanishes on $(Z\setminus X_2)\cap g^{-1}(V_1)$
but does not vanish at any point of $g^{-1}(P)$.
Since $Z\cap X_1$ is affine,  $f_P\circ g$ can be extended
to a regular function $F_p$ on $X_2$.

Set $V_P:=V_1\setminus (f_P=0)$ and
$X_P:=X_2\setminus (F_P=0)$. The restriction
$(Z\cap X_P)\to V_P$ is  finite since, by construction, 
 $X_P\cap Z$ is the preimage of $V_P$. \qed
\medskip

\begin{defn} We say that an algebraic space $X$ has the
{\it Chevalley-Kleiman property} if $X$ is separated and 
every finite subscheme is contained in an open affine subscheme.
In particular, $X$ is  necessarily a scheme.
\end{defn} 

These methods give the following interesting corollary.

\begin{cor}\label{CK.prop.desc}
 Let $\pi:X\to Y$ be a finite and surjective morphism
of separated, excellent  algebraic spaces.  
Then  $X$ has the
 Chevalley-Kleiman property iff $Y$ has.
\end{cor}

Proof. 
Assume that $Y$ has the
 Chevalley-Kleiman property and let $P\subset X$ be a finite subset.
Since $\pi(P)\subset Y$ is finite, there is an open
affine subset  $Y_P\subset Y$ containing  $\pi(P)$. Then
$g^{-1}(Y_P)\subset X$ is an open
affine subset containing  $P$.

Conversely, assume that $X$ has the
 Chevalley-Kleiman property. By the already established direction,
we may assume that $X$ is normal.
Next let $ Y^n$ be the normalization of $Y$.
Then $X\to  Y^n$ is finite and dominant.
Fix irreducible components $X_1\subset X$
and $Y_1\subset Y^n$ such that the
induced map $X_1\to Y_1$ is finite and dominant.
Let $\pi'_1:X'_1\to X_1\to Y_1$ be the Galois closure of $X_1/Y_1$
with Galois group $G$.
We already know that $X'_1$ has the
 Chevalley-Kleiman property, hence there is an
open
affine subset $X'_P\subset X_1$ containing  $\bigl(\pi'_1\bigr)^{-1}(P)$.
Then $U'_P:=\cap_{g\in G}g(X_P)\subset  X'_1$ is  affine,
 Galois invariant and 
$\bigl(\pi'_1\bigr)^{-1}\bigl(\pi'_1(U'_P)\bigr)=U'_P$.

Thus $U'_P\to \pi'_1(U'_P)$ is finite and,
by Chevalley's theorem  \cite[Exrc.III.4.2]{hartsh},
 $\pi'_1(U'_P)\subset  Y_1$
is an open
affine subset containing  $P$. 
Thus $Y^n$  has the  Chevalley-Kleiman property.

Next consider the normalization map  
 $g:Y^n\to \red Y$.
There are lower dimensional closed subschemes
$P\subset V\subset \red Y$ and 
$Z:=g^{-1}(V)\subset Y^n$ such that
$g:Y^n\setminus Z\cong \red Y\setminus V$ is an isomorphism.
By induction on the dimension, $V$ has the  Chevalley-Kleiman property.

By (\ref{affine.red.lem}) there are open  affine subsets
$P\subset V_P\subset V$ and $g^{-1}(P)\subset Y^n_P\subset Y^n$
such that $g$ restricts to a finite morphism
$g: Z\cap Y^n_P\to V_P$. Thus, by (\ref{glue.lem.affine}),
$g\bigl(Y^n_P\bigr)\subset \red Y$ is open, affine and it contains $P$.
Thus $\red Y$ has the  Chevalley-Kleiman property.

Finally,  $\red Y\to Y$ is a homeomorphism, thus
if $U\subset \red Y$ is an affine open subset
and  $U'\subset Y$  the ``same'' open subset of $Y$
then $U'$ is also affine by
Chevalley's theorem and so $Y$  has the  Chevalley-Kleiman property.
\qed
\medskip

\begin{exmp} Let $E$ be an elliptic curve and set $S:=E\times \p^1$.
Pick a general $p\in E$ and $g:E\times \{0,1\}\to E$
be the identity on $E_0:=E\times \{0\}$ and
translation by $-p$ on $E_1:=E\times \{1\}$.
Where are the affine charts on the quotient $Y$?

If $P_i\subset E_i$ are 0-cycles then
there is an ample divisor $H$ on $S$ such that
$(H\cdot E_i)=P_i$ iff  $\o_{E_0}(P_0)=\o_{E_1}(P_1)$
under the identity map $E_0\cong E_1$.

Pick any $a,b\in E_0$ and let $a+p,b+p\in E_1$ be obtained
by translation by $p$. 
Assume next that 
$2a+b=a+p+2(b+p)$, or, equivalently, that $3p=a-b$. 
Let $H(a,b)$ be an ample divisor on $S$ such that
$H(a,b)\cap E_0=\{a,b\}$ and 
$H(a,b)\cap E_1=\{a+p,b+p\}$.
Then $U(a,b):=S\setminus H(a,b)$ is affine
and $g$ maps $E_i\cap U(a,b)$ isomorphically
onto $E\setminus\{a,b\}$ for $i=0,1$. As we vary $a,b$
(subject to $3p=a-b$) we get an affine covering
of $Y$.

Note however that the curves $H(a,b)$ do not give
Cartier divisors on $Y$. In fact, for non-torsion $p\in E$,
every line bundle on $Y$ pull back from
the nodal curve obtained from the $\p^1$ factor by gluing
the points $0$ and $1$ together.
\end{exmp}

\section*{Appendix by Claudiu Raicu}

\newcommand{\x}{\textit{\textbf{x}}}
\newcommand{\y}{\textit{\textbf{y}}}
\newcommand{\zz}{\textit{\textbf{z}}}

\begin{say}\label{raicu}
Let $A$ be a noetherian commutative ring and $X=\a_S^n$ the 
$n$-dimensional affine space over $S=\spec A$. Then 
$\o_X\simeq A[\x]$, where $\x=(x_1,\dots,x_n)$.
 To give a finite equivalence relation $R\subset X\times_S X$ is 
equivalent to giving an ideal $I(\x,\y)\subset A[\x,\y]$ which 
satisfies the following properties:
\begin{enumerate}
  \item (reflexivity) $I(\x,\y)\subset (x_1-y_1,\dots,x_n-y_n)$.
  \item (symmetry) $I(\x,\y)=I(\y,\x).$
  \item (transitivity) $I(\x,\zz)\subset I(\x,\y)+I(\y,\zz)$ in $A[\x,\y,\zz]$.
  \item (finiteness) $A[\x,\y]/I(\x,\y)$ is finite over $A[\x]$.
\end{enumerate}

Suppose now that we have an ideal $I(\x,\y)$ satisfying (1-4) and 
let $R$ be the equivalence relation it defines. 
If the geometric quotient exists, then by (\ref{cat.geom.quot}) it
is of the form  $\spec A[f_1,\dots,f_m]$ for some polynomials 
$f_1,\dots,f_m\in A[\x]$. It follows that
\[I(\x,\y)\supset (f_i(\x)-f_i(\y):i=1,2,\dots,m)\]
and $R$ is {\it effective} iff equality holds.

We are mainly  interested in the case when 
$A$ is $\z$ or some field $k$ and $I$ is homogeneous. Write $I$ as
$I(\x,\y)=\bigl(J(\x,\y),f(\x,\y)\bigr)$, where $J$ is an ideal of the form
\[J(\x,\y)=\bigl(f_i(\x)-f_i(\y):i=1,2,\dots,m\bigr),\]
with homogeneous $f_i\in A[\x]$ such that 
$A[\x]$ is a finite module over $A[f_1,\dots,f_m]$ 
 and $f\in A[\x,\y]$ a homogeneous polynomial
 that satisfies the cocycle condition
$$
f(\x,\y)+f(\y,\zz)-f(\x,\zz)\in J(\x,\y)+J(\y,\zz)\subset A[\x,\y,\zz].
\eqno{(\ref{raicu}.1)}
$$
The reason we call (\ref{raicu}.1) a cocycle condition is the following. 
If we let $B=A[f_1,\dots,f_m]$, $C=A[\x]$ and consider the complex 
(starting in degree zero)
$$
C\to C\otimes_B C\to\cdots\to C^{\otimes_B m}\to\cdots
\eqno{(\ref{raicu}.2)}
$$
with differentials given by the formula
\[d_{m-1}(c_1\otimes c_2\otimes\cdots\otimes c_m)=
\sum_{i=1}^{m+1}(-1)^i c_1\otimes\cdots\otimes c_{i-1}\otimes 1\otimes c_i\otimes\cdots\otimes c_m,\]
then $C\otimes_B C\simeq A[\x,\y]/J(\x,\y)$, 
$C\otimes_B C\otimes_B C\simeq A[\x,\y,\zz]/(J(\x,\y)+J(\y,\zz))$, and if the polynomial 
$f(\x,\y)$ satisfies (\ref{raicu}.1),
 then its class in $C\otimes_B C$ is a 1-cocycle
 in the complex (\ref{raicu}.2).

Any ideal $I(\x,\y)$ defined as above is the ideal of an equivalence relation. 
To show that the equivalence relation it defines is noneffective it suffices
 to check that $f(\x,\y)$ is not congruent to a difference modulo $J(\x,\y)$. 
This can be done using a computer algebra system by computing the finite 
$A$-module $U$ of homogeneous forms of the same degree as $f$ which are
 congruent to differences modulo $J$, and checking that $f$ is not contained
 in $U$. We used Macaulay 2 to check that the following example gives a 
noneffective equivalence relation (we took $A=\z$ and $n=2$):

\[f_1(\x)=\x_1^2,\ f_2(\x)=\x_1\x_2-\x_2^2,\ f_3(\x)=\x_2^3,\]
\[f(\x,\y)=(\x_1\y_2-\x_2\y_1)\y_2^3.\]
\[I(\x,\y)=(\x_1^2-\y_1^2,\x_1\x_2-\x_2^2-\y_1\y_2+\y_2^2,
\x_2^3-\y_2^3,(\x_1\y_2-\x_2\y_1)\y_2^3).\]

We also claim that this example remains noneffective 
 after any base change. Indeed, 
the $A$-module $V$ 
generated by the forms of degree $5(=\deg(f))$ in $I$ and the differences
 $g(\x)-g(\y)$ with $g$ homogeneous of degree $5$, is a direct summand in $U$.
 Elements of $V$ correspond to 0-coboundaries in (\ref{raicu}.2). The module 
$W$ consisting of elements of $k[\x,\y]_5$ whose classes in $k[\x,\y]/J$ are
 1-cocycles is also a direct summand in $U$. The quotient $W/V$ is a free
 $\z$-module $H$ generated by the class of $f(\x,\y)$. This shows that
 $W=V\oplus H$, hence for any field $k$ we have $W_k=V_k\oplus H_k$, where 
for an abelian group $G$ we let $G_k=G\otimes_{\z}k$. If we denote by $d_i^k$
 the differentials in the complex obtained from (\ref{raicu}.2) by base 
changing from $\z$ to $k$, then we get that $\im d_0^k=V_k$ and
 $\ker d_1^k\supset W_k$. It follows that the nonzero elements of $H_k$
 will represent nonzero cohomology classes in (\ref{raicu}.2) for any field
 $k$, hence our example is indeed universal.
\end{say}

One can see \cite[Lem.4.3]{CR} that 
all  homogeneous noneffective equivalence relations
are contained in a   homogeneous noneffective equivalence relation
constructed as above.

In the positive direction, we have the following result in the toric case,
where a  \textit{toric equivalence relation} (over a field $k$)
is a scheme theoretic equivalence relation $R$ on a (not necessarily normal)
toric variety $X/k$
that  is invariant under the diagonal action of the torus.

\begin{thm}\label{thmamit} \cite[Thm.4.2]{CR} Let $k$ be a field, 
$X/k$ an  affine toric variety, and $R$ a 
toric equivalence relation on $X$. Then there exists an affine toric variety 
$Y/k$ together with a toric map $X\to Y$ such that $R\simeq X\times_Y X$.
\end{thm}

Notice that we do not require the equivalence relation to be finite.

 \begin{ack} We thank  D.~Eisenbud,  M.~Hashimoto,
 C.~Huneke, K.~Kurano, M.~Lieblich, 
P.~Roberts,  Ch.~Rotthaus, D.~Rydh  and R.~Skjelnes
for useful comments, corrections and references.
Partial financial support  was provided by  the NSF under grant number 
DMS-0758275.
\end{ack}

\bibliography{refs}

\vskip 0.5cm

\noindent Princeton University, Princeton NJ 08544

\begin{verbatim}kollar@math.princeton.edu
\end{verbatim}

\noindent University of California, Berkeley, CA 94720-3840

\noindent  Institute of Mathematics ``Simion Stoilow" of the Romanian Academy

\begin{verbatim}claudiu@math.berkeley.edu
\end{verbatim}

\end{document}